\documentclass{article}
\usepackage{amsmath}
\numberwithin{equation}{section}
\usepackage{amssymb}
\usepackage[dvipsnames]{xcolor}
\usepackage{comment}
\usepackage{hyperref}
\usepackage{graphicx}
\usepackage{yhmath}
\usepackage[margin=2.5cm]{geometry}
\usepackage{subcaption}
\newtheorem{proposition}{Proposition}
\newtheorem{definition}{Definition}

\newtheorem{theorem}{Theorem}
\newtheorem{lemma}{Lemma}
\numberwithin{equation}{section}
\numberwithin{proposition}{section}
\numberwithin{definition}{section}
\numberwithin{theorem}{section}
\numberwithin{lemma}{section}

\begin{document}

\title{A Simple Manifold with No Convex Foliation in Three Dimensions and Above}

\author{Selim Amar}

\maketitle
\noindent
Email: selama@stanford.edu \\
Department of Mathematics, Stanford University

\vspace{10pt}

\begin{abstract}

\noindent We present examples of simple metrics that admit no convex foliation with respect to the geodesic flow in 3-dimensions and above. They are constructed as oscillating perturbations of the round metric in a ball. This answers a question of A. Vasy and G. Uhlmann pertaining to their method of solving the X-ray transform. The example was found using Rethlas, a harness for AI agents, running on ChatGPT-Astra 6.

\end{abstract}

\tableofcontents
\section{Introduction}
In \cite{uhlmann2012inverseproblemlocalgeodesic}, A. Vasy and G. Uhlmann proposed a technique to invert the X-ray transform on a compact Riemannian manifold with boundary $M$ based on a geometric condition: the presence of a convex foliation. As a recall, the X-ray transform maps functions on the manifold $M$ to functions on its sphere bundle $SM$ by integrating along geodesics:
\begin{equation}
    I(f)(\beta) = \int_{\gamma_{\beta}} f \quad \beta \in SM, \frac{d}{dt}\gamma|_{t = 0} = \beta
\end{equation}
Assuming the existence of a convex foliation, which we will shortly define, they showed that one could invert the following operator
\begin{equation}
    N = L \circ \chi \circ I \quad L(g)(x) = \int_{S_x M} g d\nu \quad \chi \in C^{\infty}(SM)
\end{equation}
which gives a reconstruction mechanism for $f$ given $I(f)$. This approach based on convexity was then used in numerous works to study various inverse problems related to the X-ray transforms, such as the tensor X-ray transform in \cite{StefanovUhlmannVasy2018}, a vector-valued weighted X-ray transform in \cite{PaternainSaloUhlmannZhou2019} and the nonlinear boundary rigidity problem in \cite{StefanovUhlmannVasy2021}. The whole construction was also recast in a semiclassical framework in \cite{Vasy2024}. The crucial input was that a convex foliation effectively allowed the problem to be localized. For the purpose of this paper, a convex foliation will be given by the following definition, which may be relaxed.
\begin{definition}
    A convex foliation of $M$ is given by a smooth function $\rho: M \rightarrow [0, \infty)$ such that for some $T>0$ we have $\rho^{-1}(0) = \partial M$, $d\rho \neq 0$ for $\rho < T$ and its level sets are convex:
    \begin{equation}
        Hess (\rho)|_{ker d\rho} < 0 \quad \text{ for } \rho < T
    \end{equation}
    Furthermore, $\rho$ is trivial for $\rho \geq T$, i.e. $int(\{ \rho \geq T \}) = \emptyset$.
\end{definition}
Since the X-ray transform is relevant to the boundary rigidity problem and Michel's conjecture, which states that simple manifolds are boundary rigid, they posed a natural question: do all simple manifolds \footnote{The definition of simple manifolds we will follow is definition 3.8.1 in \cite{paternain2023geometric}} admit a convex foliation ? The purpose of this paper is to exhibit a counterexample to that question, i.e. a simple manifold with no convex foliations, in dimension 3 and higher. 
\begin{theorem}
    In dimension 3 and above, there exists explicit metrics on the ball which are simple but admit no convex foliations, of the form
    \begin{equation}
        g_{\epsilon, N} = dr^2 + \sin^2(r) h + \epsilon^2 \psi\left(\frac{\pi/2 - r}{\epsilon}\right) Q_N \quad 0 \leq r \leq \pi/2 - \epsilon
    \end{equation}
    where $Q_N$ is an oscillating tensor on the sphere and $\psi$ a smooth function. 
\end{theorem}
Two-dimensional simple manifolds are already known to possess convex functions \footnote{These are essentially convex foliation but $\rho$ is convex in all directions.}, see for example \cite{PaternainSaloUhlmannZhou2019} lemma 2.2 which is based on the curvature flow argument of \cite{BeteluGulliverLittman2002}. The example was found using Rethlas \cite{ju2026automatedconjectureresolutionformal}, a harness for AI agents, running on ChatGPT-Astra 6; this paper provides a human-readable exposition along with a simplification of the original arguments. The author's use of AI is discussed further in the corresponding section at the end of this paper. As in every other AI constructions, it can be difficult isolating where it found its inspiration from, especially after modifying the proof. It seems that the approach taken here goes back to the work of R. Gulliver in \cite{Gulliver1975}, where he showed that there exist manifolds with sectional curvature of both signs without conjugate points which still possess focal points. He did so by modifying the metric in a compact geodesic ball within the manifold:
\begin{equation}
    g = dr^2 + f(r) h \quad f(r) \text{ being the function modified}
\end{equation}
In particular, the construction interpolates between the metric within the ball and outside of it along a thin shell $(r, r+\epsilon)$ in such a way as not to create conjugate points while being strong enough to create focal points, which is exactly what the perturbed shell construction of this paper does to obstruct the existence of a convex foliation. The main difference between the two is that in our case, $Q_N$ and not $f(r)$ is what we tune to make the perturbation work, while $\psi(r)$ plays an auxiliary role. A quantitative application of R. Gulliver's method can also be found in \cite{Dibble2017} based on the same techniques, where the author showed that the convexity radius of a manifold could be made arbitrarily small compared to the injectivity radius \footnote{There the construction is quite different than this paper, going through large covers of the manifold instead while keeping the size of the perturbed ball fixed.}. Finally, the author would also like to recognize that such a construction was found independently by other researchers using AI tools, in particular Kelvin Lam, Qiuye Jia, Anton Petrunin and Dmitry Korshunov. Their construction use the same basic criteria of an embedded manifold obstructing convex foliations together with some criterion to conserve simplicity, but proceed differently to construct it. K. Lam embeds a minimal surface in a five dimensional ball using a conformal factor \cite{lam2026convexfoliations}. This provides an obstruction to a convex foliation, while the simplicity is proven (as in this paper) by establishing a positive index using quantitative estimates along great circles. Q. Jia embeds a geodesic network and uses it to construct a measure which when paired with hessians of functions vanishes, which contradicts convexity (it also makes use of \cite{PaternainSaloUhlmannZhou2019} to relate convex foliation and functions). As for Korshunov and Petrunin, the author is less familiar with their construction as it was pointed out to him through a mathoverflow post \cite{petrunin2026geodesic}, but from what he understands they also embed a geodesic net within a simple manifold, in both three and four dimensions. The following construction is more 'analytical' in nature, which will please or displease certain readers. 

\section{Criteria for Simplicity and the Absence of a Convex Foliation}

We now describe criteria for determining when a given metric admits no convex foliation, and when it is simple. By then showing that there exists examples of metrics satisfying both, we will have succeeded in proving that there exist simple metrics with no convex foliation. 

\subsection{Criterion for the Absence of a Convex Foliation}
Let's begin with the convex foliation condition:
\begin{proposition}
\label{prop:ConvexBarrier}
    Let $\Omega \subset int(M)$ be a smooth compact domain with nonempty interior. Suppose that at every point of its boundary $\partial \Omega$ there exists a nonzero tangent vector $v$ such that its extrinsic curvature $\Pi$ satisfies
    \begin{equation}
        \Pi(v, v) \leq 0
    \end{equation}
    Then $M$ admits no convex foliation. 
\end{proposition}
\textbf{Proof:} The idea is quite simple: any foliation would have to eventually touch the boundary $\partial \Omega$, at which point its convexity would fail. Let $\rho$ be a convex foliation for $M$. By compactness, there is a minimal value for $\rho$ in $\Omega$, corresponding to a first contact point $x$. It must be smaller than $T$ because the interior of $\Omega$ is nonempty. It must be at the boundary because $d\rho_x \neq 0$ there. Being a minimum at the boundary, tangential derivatives of $\rho$ must vanish there, hence $d\rho_x$ is normal to the surface $\partial \Omega$, pointing inward:
\begin{equation}
    \nabla \rho_x = -a n_x \quad n_x = \text{outward normalized normal}
\end{equation}
The statement then follows from comparing the Hessian of $\rho$ and the extrinsic curvature at $x$:
\begin{equation}
    Hess_{\partial \Omega}(\rho|_{\partial})(v, v) = Hess_M (\rho)(v, v) + a \Pi(v, v)
\end{equation}
Since $x$ is a minimum on $\partial \Omega$, the Hessian on the left is nonnegative, hence we get a lower bound
\begin{equation}
    Hess_M(\rho)(v,v) \geq -a \Pi(v, v) \geq 0
\end{equation}
which contradicts the convexity property since $v \in ker d\rho_x$, being tangent to $\partial \Omega$.
$\Box$
\vspace{0.2in}
\newline
The key idea is that this is a \textit{pointwise criterion}; it only needs a single direction to fail in order to obstruct convexity. In terms of a metric perturbation $Q$, it will correspond to one of its eigenvalues being large enough to induce a negative direction for $\Pi$. On the other end, making sure that the resulting metric stays simple is an \textit{integrated condition}, which allows us to use cancellation to enforce it. Note that we will prove something stronger than what is necessary to apply this criterion, as for us $\Pi(v, v) < 0$, i.e. we will construct directions of 'strict concavity', not just failure of convexity.

\subsection{Criterion for Simplicity}
The criterion for simplicity we use is a \textit{positive Dirichlet index} criterion for boundary chords \footnote{Boundary chords will be defined as unit-speed geodesics starting inward from the boundary and ending at the next boundary contact.}, which follow from proposition 4.7 and theorem 4.5 of \cite{GulliverLittman2000}. For the reader unfamiliar, the Dirichlet index is the second variation of the energy of a family of curves $F(s, t)$ with fixed endpoints:
\begin{equation}
     E''(0) = I_{\gamma}(J, J) = \int_0^T ( |D_t J|^2 - \left< \mathcal{R}J, J \right>)dt \quad E(s) = \frac{1}{2} \int_0^T |\partial_t F(s, t)|^2 dt \quad J = \frac{\partial F}{\partial s}|_{s = 0}
\end{equation}
where $\mathcal{R}$ is the Riemann curvature tensor, and $J(t)$ the variation vector field which is normal to $\dot{\gamma}$ and vanish at the endpoints $0, T$. In particular, if the endpoints $\gamma(0)$, $\gamma(T)$ are conjugate points, then there exists a Jacobi field $J$ vanishing at the endpoints and in that case $I_{\gamma}(J, J) = 0$ since one use integration by parts to write
\begin{equation}
    \begin{aligned}
    I_\gamma(J,J) =  \langle D_tJ,J\rangle|_0^T -\int_0^T \left\langle D_t^2J+R(J,\dot\gamma)\dot\gamma,J\right\rangle\,dt = 0
\end{aligned}
\end{equation}
Thus, to avoid the presence of conjugate points along the curve, a simple criterion is that $I_{\gamma}(J, J) > 0$ for any possible variation vector field $J$ (this is discussed in more details in \cite{paternain2023geometric} section 3.7). The next proposition says that a uniform version of this over the boundary also implies simplicity.
\begin{proposition}
    Let $M$ be a smooth compact connected Riemannian manifold of dimension $\geq 2$ with strictly convex boundary. Suppose that the first exit time $T$ of boundary chords is uniformly bounded, and that for every boundary chords $\gamma$ we have the following lower bound for every nonzero normal field vanishing at the endpoints:
    \begin{equation}
        I_{\gamma}(J, J) = \int_0^T ( |D_t J|^2 - \left< \mathcal{R}J, J \right>)dt > 0
    \end{equation}
    Then $M$ is simple.
\end{proposition}
\textbf{Proof}: As explained above, this is an established criterion which follow from proposition 4.7 and theorem 4.5 of \cite{GulliverLittman2000}. These hypothesis implies theirs because we assume an uniformly bounded exit time, while strict positivity of the index excludes conjugate points. Theorem 4.5 and proposition 4.7 then give for each interior geodesic a unique extension to a (nondegenerate) boundary chord, which implies that the manifold is simple by theorem 3.8.2 in \cite{paternain2023geometric}. $\Box$
\vspace{0.2in}
\newline
As we claimed above, this is an \textit{integrated condition} hence for highly oscillatory perturbations, can be arranged through cancellations. The strategy of the paper will thus be to exhibit a metric perturbation which satisfies both criteria for non-convexity and simplicity.
\section{Concrete Oscillating Example}
We now construct a concrete example of a metric on the ball satisfying both criteria. We construct it as a perturbation of the spherical metric on a ball in geodesic coordinates:
\begin{equation}
    g_{\epsilon, N} = dr^2 + H(r)_{\epsilon, N} = dr^2 + \sin^2(r) h + \epsilon^2 \psi\left(\frac{\pi/2 - r}{\epsilon}\right) Q_N \quad 0 \leq r \leq \pi/2 - \epsilon
\end{equation}
where $\psi$ is a function supported in $(3/2, 5/2)$ equal to $0$ at $2$ that restricts the perturbation near the boundary $r \sim \pi/2$ while $Q_N$ is an oscillating tensor with frequency given by some integer $N$ constructed as a symmetric gradient of some vector field on the sphere $Q_N = sym \nabla V_N$ \footnote{One does not need to construct $Q_N$ as a symmetric gradient, but it greatly simplifies checking the integrated condition since we can apply integration by parts. } and $h(X,X)$ is the metric on the sphere. Note that $H(r)$ is the induced metric on the spheres $r = const.$, hence the perturbation only differs near the shells close to the boundary. The key properties we need are the following uniform properties over the ball
\begin{equation}
\label{eq:PropertiesNecessary}
    \lambda_{\max}(Q_N) \geq a > 0 \quad \| Q_N \|_h \leq 1 \quad lim_{N \rightarrow \infty} \| V_N \|_{\infty} = 0 \quad \psi'(2) > 4/a
\end{equation}
An example of a vector field $V_N$ on the two-dimensional sphere $\mathbb{S}^2$ (hence in three dimensions) satisfying these for $a = 1/4$ is
\begin{equation}
    V_N = \frac{1}{3} \left( \frac{\sin(2N\theta)}{2N} \partial_{\theta} + 4 \sin(\theta) \int_{\pi/2}^{\theta} \frac{q_N(s)}{\sin(s)} ds \frac{\partial_{\phi}}{\sin(\theta)} \right) \quad q_N(\theta) = \chi(8N\theta - 1) \chi(8N (\pi - \theta) - 1) \sin(N\theta)
\end{equation}
where $0 \leq \chi(t) \leq 1$ is some cutoff equal to $0$ for $t \leq 0$ and $1$ for $t \geq 1$. The proof that it satisfies these properties is relegated to appendix \ref{App:ProofEstimatesV}, along with how to make it work in any dimensions $n \geq 3$. In the next subsection, we show that if $Q_N$ satisfies the properties above \ref{eq:PropertiesNecessary}, then for sufficiently large $N$ and small $\epsilon$, the metric will satisfy both criteria simultaneously.

\subsection{The Metric admits no Convex Foliations}
The easiest one to check is the absence of convex foliation. We need to compute the second fundamental form of $g_{\epsilon, N}$ on coordinate spheres $r = const.$. It is given by
\begin{equation}
    \Pi_r = \frac{1}{2} H_{\epsilon, N}'(r) = \sin(r) \cos(r) h - \frac{\epsilon}{2} \psi' Q_N
\end{equation}
At $R = \pi/2 - 2\epsilon$, pick an eigenvector $v$ of $Q_N$ with eigenvalue $\lambda_{\max} \geq a$, so that $Q_N(v, v) \geq a$. Then 
\begin{equation}
    \Pi_R(v, v) = \sin(R) \cos(R) - \frac{\epsilon}{2} \psi'(2) \lambda_{\max} \leq \epsilon \left(2 - \frac{\psi'(2) a}{2} \right) < 0
\end{equation}
Thus, the ball with boundary $R$ provides the embedded barrier to a convex foliation from proposition \ref{prop:ConvexBarrier}. 
\subsection{The Metric is Simple}
\begin{figure}[h]
    \centering
    \includegraphics[width=0.85\linewidth]{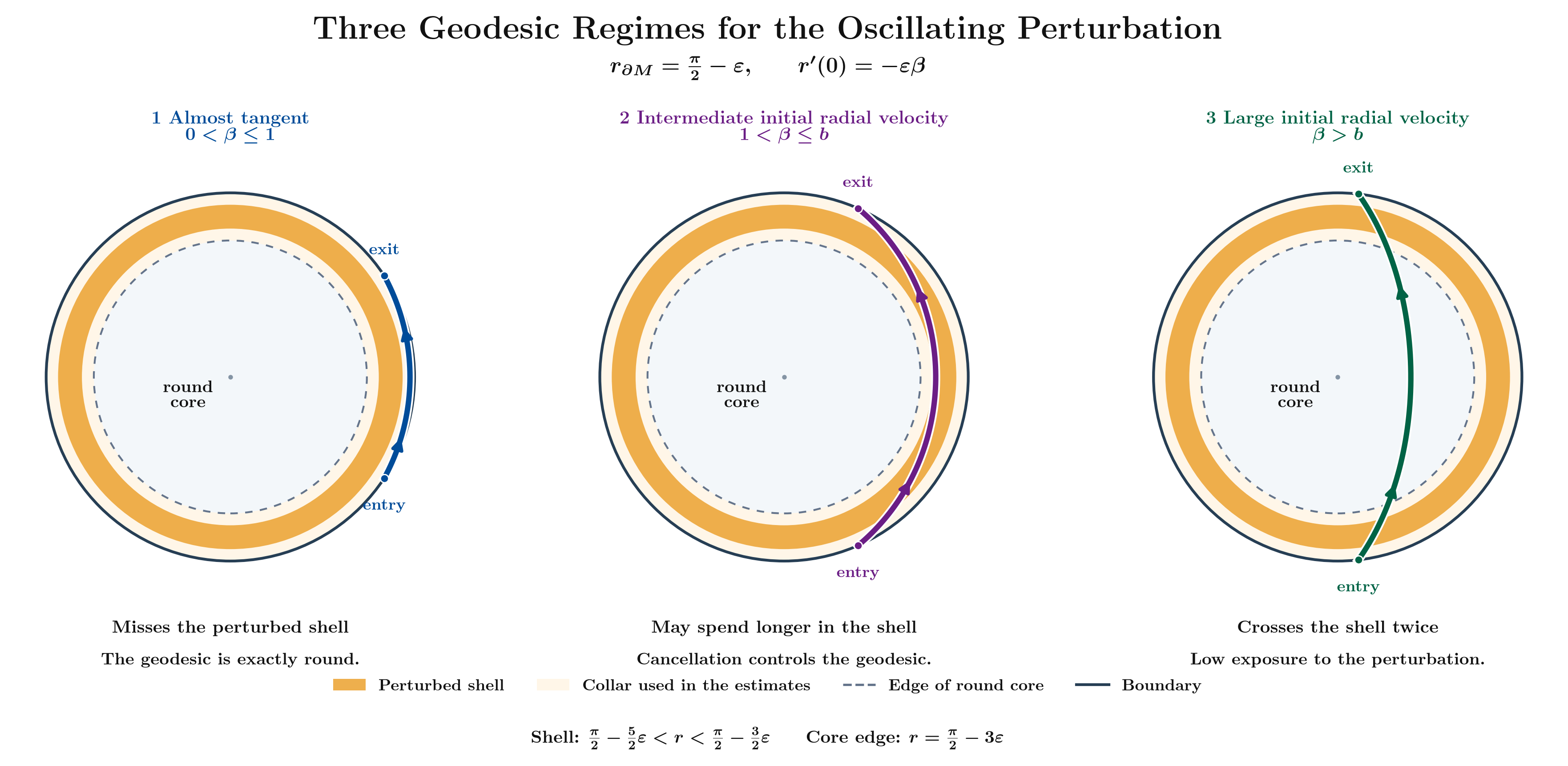}
    \caption{An illustration of the three regimes for the boundary chords. Each require a separate analysis to show that it has positive index.}
    \label{fig:ThreeGeodesics}
\end{figure}
The more involved computation is to check that these properties of $Q_N$ make the metric simple. This is where we use its norm estimates, and that it is a symmetric gradient. We will show that if the perturbation $Q$ has a cancellation property, then the metric is simple. It amounts to the following:
\begin{proposition}
    Let $Q = Q_{N}$ be a metric perturbation as above with $\| Q \|_h \leq 1$. Then there is a $\delta>0$ depending only on the cutoff $\psi$ and the dimension such that if for all unit-speed $h-$geodesics $\alpha$ on the sphere we have
    \begin{equation}
         \quad \sup_{\alpha, 0 \leq L \leq \pi} \left| \int^L_0 Q(\dot{\alpha}, \dot{\alpha}) dt \right| + \epsilon \max_{j \leq 2} \| (\nabla^h)^j Q \|_{\infty} + \epsilon < \delta
    \end{equation}
    then all boundary chords have uniformly bounded exit time and positive index. Thus, $M$ is simple. 
\end{proposition}
Our perturbation $Q_N$ automatically satisfies this for large enough $N$ since
\begin{equation}
    \sup_{\alpha, 0 \leq L \leq \pi} \left| \int^L_0 Q(\dot{\alpha}, \dot{\alpha}) dt \right| = \sup_{\alpha, 0 \leq L \leq \pi} \left| \int^L_0 \frac{d}{dt} \left<V, \dot{\alpha} \right> dt \right| \leq 2 \| V_N \|_{\infty} \rightarrow 0     
\end{equation}
hence by taking $\epsilon$ small enough (depending on $N$, which is taken large) we can arrange it to satisfy the desired bound. The proof separates the boundary chords into three cases depending on their initial (rescaled) radial velocity $\beta = -\frac{\dot{r}(0)}{\epsilon}$, $r(0) = \pi/2 - \epsilon$:
\begin{itemize}
    \item The boundary chord is almost tangent to the spheres: $ 0 \leq \beta \leq 1$.
    \item The boundary chord may spend a long time near the perturbed shell: $ 1 < \beta \leq b$.
    \item The boundary chord crosses the shell steeply: $\beta > b$.
\end{itemize}
The main idea is to compare the resulting unit-speed geodesics $\gamma(t) = (r(t), x(t))$ to round geodesics in three different ways. For the first case, the geodesic never reach the perturbed shell hence it will be a geodesic with respect to the unperturbed round metric $dr^2 + \sin^2(r) h$. It then trivially satisfies the required exit time and index properties. The second case is the most challenging; the geodesic will suffer the full effect of the perturbation. This is where we make use of the cancellation properties of our perturbation $Q$, to say that the overall \textit{integrated} effect is nonetheless perturbative relative to the round geodesics of the spherical metric. Finally, the last case is also quite simple, since by taking the initial radial speed large, the geodesic simply crosses the perturbed region rapidly hence make the perturbation harmless. In reading the proof, it is worth keeping in mind two facts about the unperturbed round metric. The first, is that any boundary chords exit in time at most $\pi$ from the sphere, which follows from the exit time being given by \footnote{This formula will come out in the proof as well, when considering the height of geodesics $y(t) = \cos(r(t))$.}
\begin{equation}
    l = 2 \arctan(\epsilon \beta \cot(\epsilon)) \leq 2 \arctan(\beta) < \pi
\end{equation}
The second is that the index has an explicit positive lower bound. Indeed, the sectional curvature is 1 hence
\begin{equation}
    I_{\gamma}(J, J) = \int_0^l |D_tJ|^2 - |J|^2 dt \geq \left(1 - \frac{l^2}{\pi^2} \right) \int_0^l |D_t J|^2 dt > 0
\end{equation}
where we've made use of Poincaré inequality
\begin{equation}
    \int_0^l |J|^2 dt \leq \frac{l^2}{\pi^2} \int_0^l |D_t J|^2 dt \quad J(0) = J(l) = 0
\end{equation}
The key is that our perturbation will preserve these up to small errors.
\vspace{0.1in}
\newline
\textbf{Proof:}
\newline
Step 1) $0 \leq \beta \leq 1$: These geodesics turn out to be the usual geodesics in the (unperturbed) ball and do not interact with the shell. Indeed, their height satisfies the ODE $y''+y = 0$ with $y(0) = \sin(\epsilon)$ and $y'(0) = \epsilon \beta \cos(\epsilon)$ hence $y = \cos(r(t))$ is
\begin{equation}
    \cos(r(t)) = \sin(\epsilon) \cos(t) + \epsilon \beta \cos(\epsilon) \sin(t)
\end{equation}
which for $\beta \leq 1$ means that $r(t)$ stays bounded away from the perturbed shell $r = \pi/2 - 3/2 \epsilon$:
\begin{equation}
    r(t) \geq \arccos\left(\sqrt{\sin^2(\epsilon) + \epsilon^2 \beta^2 \cos^2(\epsilon)}\right) > \frac{\pi}{2} - \frac{3}{2} \epsilon
\end{equation}
In particular, it is a round geodesic hence has exit time bounded by $\pi$ and positive index.
\vspace{0.1in}
\newline
Step 2) $1 < \beta \leq b$: This is the most difficult one. The idea of this step is again to compare the geodesics to the unperturbed geodesics. This time they will not be exactly round, but the integrated estimate on the perturbation $Q$ will keep them close, hence preserve their exit time bounded by $\pi$ and the index positive. This is the case where we make the most use of the oscillatory nature of our perturbation $Q$.
\begin{lemma}
    For $1 < \beta \leq b$, the unit-speed geodesics of the perturbed metric stay uniformly close to the unperturbed ones up to their exit time, in the sense that
    \begin{equation}
        \left| r(t) - \frac{\pi}{2} \right| + |\dot{r}(t)| \leq C_b \epsilon \quad |x(t) - \alpha(t)| + \left| X(t) - \dot{\alpha}(t) \right| \leq C_b \epsilon \delta
    \end{equation}
    where $\alpha$ is the $h-$unit-speed geodesic with respect to the spherical metric $h$.
\end{lemma}
We start by showing that $r$ stays roughly constant, with $\dot{r}$ small. Consider again the height $y(t) = \cos(r(t))$. Using that the radial geodesic equation is
\begin{equation}
    r''(t) = \frac{1}{2} H'(X, X) = \sin(r) \cos(r) h(X, X) - \frac{\epsilon}{2} \psi'\left( \frac{\pi/2 - r}{\epsilon} \right) Q(X, X)
\end{equation}
we derive that $y(t)$ satisfies a perturbed harmonic equation
\begin{equation}
    y'' + y = \left( \epsilon^2 \cos(r) \psi\left( \frac{\pi/2 - r(t)}{\epsilon} \right) + \frac{\epsilon}{2} \sin(r) \psi'\left( \frac{\pi/2 - r(t)}{\epsilon} \right) \right) Q(X, X)
\end{equation}
In particular the forcing is $\lesssim \epsilon$ (this uses our bound $\|Q \| \leq 1$), so pairing this with $y(0) = \sin(\epsilon)$ and $y'(0) = \epsilon \beta \cos(\epsilon)$ gives a $b-$dependent bound
\begin{equation}
    |y| + |y'| \leq C_b \epsilon \implies \left| r(t) - \frac{\pi}{2} \right| + |\dot{r}(t)| \leq C_b \epsilon
\end{equation}
We next compare the angular motion $x(t)$. The angular geodesic equation for the angular velocity $X(t)$ is
\begin{equation}
    \nabla^h_{\dot{\gamma}(t)} X = -H^{-1} (\partial_r H) r' X - (\nabla^{H(r)}_X X - \nabla^h_X X)
\end{equation}
where $\nabla$ denote Levi-Civita connections associated to the corresponding metrics. Now by what we just showed for $r(t)$, the induced metric on the spheres satisfy
\begin{equation}
    H = h + O_b(\epsilon^2) \quad \partial_r H = O_b(\epsilon)
\end{equation}
This implies that the right hand-side of this angular geodesic equation is $O_b(\epsilon \delta)$ (note that we use the control on $\nabla^h Q$ we assumed to control the difference of connections):
\begin{equation}
    |\nabla^h_{\dot{\gamma}(t)} X| \leq C_b \epsilon \delta
\end{equation}
Since unit-speed with respect to $dr^2 + H$ implies the same up to $\epsilon^2$ for $h$, i.e. $h(X, X) = 1 + O_b(\epsilon^2)$, we conclude that if $\alpha(t)$ denotes a unit-speed round angular geodesic with $\alpha(0) = x(0)$, $\dot{\alpha}(0) = X(0) / |X(0)|_h$, its difference with $(x(t), X(t))$ satisfies a forced harmonic equation with forcing bounded by $\lesssim \epsilon \delta$, hence
\begin{equation}
    |x(t) - \alpha(t)| + \left| X(t) - \dot{\alpha}(t) \right| \leq C_b \epsilon \delta
\end{equation}
$\Box$
\newline
Thus, both the radial and angular motion are uniformly perturbative. This allows us to now control the effect of $Q$ along the perturbed geodesics. Define the integrated contribution along the geodesic by
\begin{equation}
    P(t) = \int_0^t Q_{x(s)}(X(s), X(s)) ds
\end{equation}
By what we showed and our control on the derivatives of $Q$, $|Q_x(X, X) - Q_{\alpha}(\dot{\alpha}, \dot{\alpha})| \leq C_b \delta^2$ hence this contribution is comparable to the one on $h-$geodesics. Thus, by our assumptions on the integral over the $h-$geodesics, we have
\begin{equation}
    |P(t)|_{\infty} \leq \delta + C_b \delta^2
\end{equation}
This is the key estimate to control the exit time and index of the geodesics. In the following, we will use a weighted version of this estimate; for a differentiable weight $w$, integration by parts also gives
\begin{equation}
    \left| \int_0^t w(s) Q_{x(s)}(X(s), X(s)) ds \right| \leq (\delta + C_b \delta^2) \left( |w(t)| + \int_0^t |w'|ds \right)
\end{equation}
We start by showing that the exit time is still below $\pi$. Again, consider the equation for the height $y(t)$. Denote the equations by
\begin{equation}
    y'' + y = \epsilon a(t) Q_{x(t)}(X(t), X(t)) \quad a(t) = \frac{1}{2} \sin(r(t)) \psi' + \epsilon \cos(r(t)) \psi
\end{equation}
If we denote by $y_0(t) = \sin(\epsilon) \cos(t) + \epsilon \beta \cos(\epsilon) \sin(t)$ the unperturbed geodesic, then the difference of the two satisfies
\begin{equation}
    y(t) - y_0(t) = \epsilon \int^t_0 \sin(t-s) a(s) Q_{x(s)}(X(s), X(s)) ds
\end{equation}
Using that $|a| + |a'| \leq C_b$, we can apply our weighted estimate for $P(t)$ to control this difference and its derivative
\begin{equation}
    \| y(t) - y_0(t) \|_{C^1} \leq C_b \epsilon (\delta + C_b \delta^2)
\end{equation}
This implies that the geodesic $\gamma(t)$ leaves the ball before time $\pi$ at 
\begin{equation}
    l_{\ast} = \frac{\pi}{2} + \arctan(b) < \pi
\end{equation}
since at $t = l_{\ast}$, the unperturbed geodesic is outside the ball (whose boundary is $y = \sin(\epsilon)$) by
\begin{equation}
    y_0(l_{\ast}) = \frac{-b \sin(\epsilon) + \epsilon \beta \cos(\epsilon)}{\sqrt{1 + b^2}} \leq 0
\end{equation}
hence by our bound on their difference with the perturbed geodesic we have
\begin{equation}
    y(l_{\ast}) - y_0(l_{\ast}) \leq C_b \epsilon(\delta + C_b \delta^2)
\end{equation}
which for small $\delta$ implies $y(l_{\ast}) < \sin(\epsilon)$. As for the index being positive, this relies on estimating the curvature tensor. We delegate the computations to appendix \ref{App:ProofCurvature} and here simply quote the end result.
\begin{lemma}
\label{lemma:CurvatureTensor}
    Let $V$ be a parallel vector field along the geodesic $\gamma$ whose intial value is normal and almost radial \footnote{Introducing this $V$ is optional. One could instead work with $V = \partial_r$, with the caveat that the integration by parts done to estimate the oscillating part $\sim Q$ now comes with an error $O(\epsilon)$.}
    \begin{equation}
        V(0) = \frac{\partial_r - \dot{r}(0) \dot{\gamma}(0)}{\sqrt{1 - \dot{r}(0)^2}}
    \end{equation}
    Then the curvature tensor along normal vector fields is given by
    \begin{equation}
        \langle \mathcal{R}(J, \dot{\gamma}) \dot{\gamma}, J \rangle = |J|^2 - \frac{1}{2} \psi''\left( \frac{\pi/2 - r}{\epsilon} \right) Q_{x(t)}(X, X) \langle J, V \rangle^2 + \langle\mathcal{E}J, J\rangle
    \end{equation}
    where $|\mathcal{E}| \leq C_b \delta$. In other words, it is the unperturbed curvature tensor, plus an oscillating part and a pointwise small contribution.
\end{lemma}
From this formula, we have
\begin{equation}
    I_{\gamma}(J, J) = \int_0^l |D_t J|^2 dt - \int_0^l |J|^2 dt + \frac{1}{2} \int_0^l \psi'' Q_{x(t)}(X, X) \langle J, V \rangle^2 dt - \int^l_0 \langle EJ, J \rangle dt
\end{equation}
where $l$ is the length of the geodesic, which is uniformly bounded by $l_{\ast} < \pi$ by what we showed. By Poincaré inequality, we have 
\begin{equation}
    \int_0^l |J|^2 dt \leq \frac{l^2}{\pi^2} \int_0^l |D_t J|^2 dt
\end{equation}
hence we get a lower bound for the first two terms (which correspond to the unperturbed metric) by
\begin{equation}
    \int_0^l |D_t J|^2 dt - \int_0^l |J|^2 dt \geq \left(1 - \frac{l^2}{\pi^2} \right) \int^l_0 |D_t J|^2 dt
\end{equation}
As for the remaining two error terms, the one involving $Q$ can be bounded again with our integrated weighted estimate plus the fact that $J$ vanish at the endpoints:
\begin{equation}
\begin{aligned}
    \left|\int_0^l \psi'' Q_{x(t)}(X, X) \langle J, V \rangle^2 dt\right| &\leq (\delta + C_b \delta^2) \left(\int^l_0 2 |\psi'' \langle D_t J, V \rangle \langle J, V \rangle| dt + \int_0^l \left| \frac{\dot{r}}{\epsilon} \psi''' \langle J, V \rangle^2 \right| dt \right) \\
    &\leq C_b (\delta + C_b \delta^2) \left( \int_0^l |J|^2 dt + \int_0^l |J| |D_t J| dt \right) \\
    &\leq C_b (\delta + C_b \delta^2) \int_0^l |D_t J|^2 dt
\end{aligned}
\end{equation}
while the one involving $\mathcal{E}$ is already bounded by $C_b \delta$. Thus, we overall get
\begin{equation}
    I_{\gamma}(J, J) \geq \left(1 - \frac{l_{\ast}^2}{\pi^2}  - C_b(\delta + C_b \delta^2) \right) \int^l_0 |D_t J|^2 dt
\end{equation}
which for sufficiently small $\delta$ is positive.
\vspace{0.1in}
\newline
Step 3) $\beta > b$: For this case, the strategy will be to show that the geodesic spends little time in the perturbed region near the boundary $\pi/2 - 3\epsilon \leq r \leq \pi/2 - \epsilon$. It goes through three stages: first entering the perturbed shell, then going through the unperturbed center of the ball, and finally re-entering the shell on its exit. Because the speed is so large, the effect of the perturbed region will be small, which again preserves the exit time and positive index. We start by showing that the exit time, which is the time spent in each phase, is bounded by $\pi - O(\beta^{-1})$.
\begin{lemma}
    Let $l_1$, $l_2$ and $l_3$ be the time in each phases respectively. We have
    \begin{equation}
        l_1 + l_3 \leq \frac{4}{\beta} + O(\beta^{-3}) \quad l_2 \leq \pi - \frac{6}{\beta} + O(\beta^{-3})
    \end{equation}
    In particular the total time $l_1 + l_2 + l_3$ is at most $\pi - \frac{1}{\beta}$ for $\beta$ large.
\end{lemma}
The bound on $l_1, l_3$ is a simple derivative estimate. The radial geodesic equation there satisfied $|r''| \leq C\epsilon$ over a width $2\epsilon$ region. Since initially $(r')^2 = \epsilon^2 \beta^2$, we have
\begin{equation}
    |\frac{d}{dr}(r'^2)| = |2r''| \leq C \epsilon \implies (r')^2 \geq \epsilon^2(\beta^2 - M)
\end{equation}
which for large $\beta$ never crosses $0$. Thus, the geodesics goes through the shell, and in time
\begin{equation}
    l_1, l_3 \leq \frac{2 \epsilon}{\epsilon \sqrt{\beta^2 - 2M}} = \frac{2}{\beta} + O(\beta^{-3})
\end{equation}
As for the time spent during the core, it is then a unperturbed geodesic with inward radial speed 
\begin{equation}
    \epsilon^2 (\beta^2 - M) \leq w^2 \leq \epsilon^2(\beta^2 + M)
\end{equation}
Its height $y(t)$ then satisfies
\begin{equation}
    y(t) = \cos(r(t)) = \sin(3\epsilon) \cos(t) + w \cos(3\epsilon) \sin(t)
\end{equation}
hence it returns to the perturbed shell at $r = \pi/2 - 3\epsilon$ at
\begin{equation}
    2 \arctan(w \cot(3\epsilon)) \leq 2 \arctan\left( \frac{\sqrt{\beta^2 + M}}{3} \right) = \pi - \frac{6}{\beta} + O(\beta^{-3})
\end{equation}
This shows that the geodesics exit in bounded time $\Box$. 
\vspace{0.1in}
\newline
As for the index staying positive, by the computations we did for step 2) we have
\begin{equation}
    I_{\gamma}(J, J) \geq \left( 1 - \frac{l^2}{\pi^2} \right) \int^l_0 |D_t J|^2 dt - (\kappa - 1) \int_{shell} |J|^2 dt
\end{equation}
where $\kappa \geq 1$ bounds the curvature in the shell \footnote{Note that this bound is independent of $b$, and only depends on $\psi$, $\delta$ and the dimension through estimating $H$ and its derivatives.}. Since $J$ vanishes at the endpoints, we have
\begin{equation}
    |J(t)|^2 \leq t \int_0^l |D_t J|^2 dt \quad |J(t)|^2 \leq (l - t) \int_0^l |D_t J|^2 dt
\end{equation}
hence the integral over the shell simplifies to
\begin{equation}
    \int_{shell} |J|^2 dt \leq \frac{l_1^2 + l_3^2}{2} \int_0^l |D_t J|^2 dt 
\end{equation}
The total lower bound is thus
\begin{equation}
    I_{\gamma}(J, J) \geq \left( 1 - \frac{l^2}{\pi^2} - \frac{(\kappa - 1)(l_1^2 + l_3^2)}{2} \right) \int_0^l |D_t J|^2 dt \geq \frac{c}{\beta} \int_0^l |D_t J|^2 dt > 0
\end{equation}
This concludes the proof, by choosing $b$ large enough so that step 3) holds (i.e. large $\beta$), and then choosing $\delta$ small enough so that step 1) and 2) hold. $\Box$.
\vspace{0.2in}
\newline
We finish with some remarks about the proof. First, the obstruction holds for any dimension, as long as the perturbation $Q_N$ satisfies the required estimates. It can thus obstruct attempts at using convexity to study the X-ray transform on general simple manifolds, by constructing perturbations which are concave in multiple dimensions on some embedded sphere. Indeed, one can construct $Q_N$ still satisfying the desired estimates and with multiple eigenvalues having lower bound $\lambda \geq a$. Such metrics would have multiple directions that are concave, which makes it difficult to use convexity, even in higher dimensions. The extreme case of this is a 'codimension 2' concave subspace at every point (which is the maximal amount of concave dimensions that such a $Q_N$ can admit), which would obstruct the standard convex-foliation approach for simple manifolds in any dimensions greater than $3$, even after throwing away some of the overdetermined information. Such a $Q_N$ can be constructed, and is discussed rapidly in appendix \ref{App:ProofEstimatesV}. Secondly, we'd like to conclude this text with the key takeaway from this construction, which will be familiar to people who've taught about this problem before: \textbf{convexity is a pointwise condition, while simplicity is an integrated condition}. The two can thus hold in separate instances when cancellation occurs.
\newpage
\section{Acknowledgments}

I would like to thank my advisor Andras and my director of graduate studies Otis for encouraging me to decipher the meaning of this construction. Furthermore, I am grateful to the anonymous donor who sponsored our department's subscriptions to ChatGPT, without which this counterexample would not have been found. 

\section{Declarations}
\textbf{Use of AI}
\vspace{0.1in}
\newline
\noindent As explained in the introduction, ChatGPT Astra-6 was used to find the counterexample and handle the technical estimates, following a discussion with the author's advisor Andras Vasy. The first trials were mostly unsuccessful. It was after the author made use of the harness Rethlas \cite{ju2026automatedconjectureresolutionformal} that the first counterexample was found. It was initially quite difficult to comprehend. The author's main contributions was working with the system to simplify the proof, in particular proposing the idea of using integration by parts to simplify the integrated estimates (which amounted to using $Q_N$ as a symmetric derivative + the weighted integrated estimate) along with extracting the main ideas found in this text (in particular the comparison with unperturbed geodesics in different regimes), which were hidden in the original AI's exposure. The figure was also done by ChatGPT Astra-6, and an initial draft of the paper was reviewed by it to catch typos. Furthermore, it was also used to survey the literature for precedents to the ideas present in this paper. This paper contains no AI-written content.
\vspace{0.2in}
\newline
\noindent \textbf{Funding}
\vspace{0.1in}
\newline
\noindent The research leading to these results received funding from the FRQ (' Fonds de Recherche du Québec') under grant number \#332587 and from the NSF ('National Science Foundation') under grant numbers DMS-2247004 and DMS-2553664.
\vspace{0.2in}
\newline
\noindent \textbf{Conflict of Interest}
\vspace{0.1in}
\newline
\noindent The author has no relevant financial or non-financial interests to disclose.
\vspace{0.2in}
\newline
\noindent \textbf{Ethical Statement}
\vspace{0.1in}
\newline
\noindent This research did not involve human or animal participants and therefore was not subjected to compliance with ethical standards concerning them.
\vspace{0.2in}
\newline
\noindent \textbf{Data Availability}
\vspace{0.1in}
\newline
\noindent The author has no relevant data to share.

\bibliographystyle{acm}
\bibliography{main}

\appendix
\section{Proof that $V_N$ satisfies the claimed properties.}
\label{App:ProofEstimatesV}
Here we present the proof that the following vector field on the sphere with metric $h = d\theta^2 + \sin(\theta)^2 d\phi^2$
\begin{equation}
\begin{aligned}
    &V_N = \frac{1}{3} \left( f_N(\theta) \partial_{\theta} + g_N(\theta) \frac{\partial_{\phi}}{\sin(\theta)} \right) \\
    & f_N(\theta) = \frac{\sin(2N\theta)}{2N} \quad g_N(\theta) =  4 \sin(\theta) \int_{\pi/2}^{\theta} \frac{q_N(s)}{\sin(s)} ds \\
    & q_N(\theta) = \chi(8N\theta - 1) \chi(8N (\pi - \theta) - 1) \sin(N\theta) 
\end{aligned}
\end{equation}
satisfies the desired property as an oscillating perturbation. Note that it is a smooth vector field, even at the poles. We start by showing that $V_N \rightarrow 0$ as $N \rightarrow \infty$. We trivially have $|f_N| \leq 1 / N$, so we're left to show $g_N$ can be bounded too. First consider the region where the cutoffs $\chi$ are $1$, say $1 / (4N) \leq \theta \leq \pi/2$. Then integrating by part the oscillating factor $\sin(Ns)$ gives
\begin{equation}
    |g_N| = \left|4 \sin(\theta) \int_{\theta}^{\pi/2} \frac{\sin(Ns)}{\sin(s)} ds \right| \leq \frac{4\sin(\theta)}{N} \left( 1 + \frac{1}{\sin(\theta)} + \int_{\theta}^{\pi/2} \frac{\cos(s)}{\sin^2(s)} ds \right) = \frac{8}{N}
\end{equation}
For $\theta \leq 1/(4N)$, we instead use that
\begin{equation}
    \frac{|q_N(s)|}{\sin(s)} \leq N
\end{equation}
which gives
\begin{equation}
    |g_N| \leq 4 \sin(\theta) \left( \frac{2}{N \sin(1/ (4N))} + N \left(\frac{1}{4N} - \theta \right) \right) \leq \frac{8}{N} + \frac{1}{4N} = \frac{33}{4N}
\end{equation}
The other angles proceed by reflection about $\pi/2$. Thus, $V_N$ goes to $0$ at rate $1/N$. Next we compute its symmetric derivative
\begin{equation}
    sym \nabla V_n = \frac{1}{3}
\begin{pmatrix} f'(\theta) & \frac12\left(g'(\theta) - g(\theta)\cot\theta\right) \\ \frac12\left(g'(\theta)-g(\theta)\cot\theta\right) & f(\theta)\cot\theta
\end{pmatrix}
\end{equation}
The off-diagonal terms simplify using that $g_N' - g_N \cot(\theta) = 4 q_N$. We thus get
\begin{equation}
    Q_N = \frac{1}{3}
\begin{pmatrix} \cos(2N\theta) & 2 q_N(\theta) \\ 2 q_N(\theta) & \frac{\sin(2N\theta) \cot(\theta)}{2N}
\end{pmatrix}
\end{equation}
We have to check that this has norm bounded by $1$, and a lower bounded eigenvalue. For the norm, note that all the coefficients are bounded by at most $1$, hence
\begin{equation}
    |Q_N((u, v), (u, v))| = \frac{1}{3} \left| \cos(2N\theta) u^2 + 4 q_N uv + \frac{\sin(2N\theta) \cot(\theta)}{2N} v^2 \right|  \leq \frac{1}{3} (u^2 + 4 |u v| + v^2) \leq u^2 + v^2
\end{equation}
As for the eigenvalue, note that its characteristic polynomial is
\begin{equation}
    p(\lambda) = \left(\lambda - \frac{\cos(2N\theta)}{3}\right) \left(\lambda - \frac{\sin(2N\theta) \cot(\theta)}{6N}\right) - \frac{4}{9} q_N^2
\end{equation}
Again in the region where the $\chi$ cutoffs are 1, the expression simplifies to (using $4q_N^2 = 2(1 - \cos(2N\theta))$)
\begin{equation}
\begin{aligned}
    p(1/3) &= \frac{1}{9} \left(\left(1 - \cos(2N\theta) \right)\left(1-\frac{\sin(2N\theta) \cot(\theta)}{2N}\right) - 2(1 - \cos(2N\theta)) \right) \\
    & = \frac{-1}{9} \left(1 - \cos(2N\theta) \right) \left( 1 + \frac{\sin(2N\theta) \cot(\theta)}{2N} \right) \leq 0
\end{aligned}
\end{equation}
Thus the maximal eigenvalue is $\lambda_{\max} \geq 1/3$. As for the other angles, say $0 \leq \theta \leq 1/4N$, in that case
\begin{equation}
    \cos(2N\theta) \geq 1 - \frac{(2N\theta)^2}{2} \geq \frac{7}{8}
\end{equation}
hence now if we evaluate the polynomial at $\lambda = \cos(2N\theta) / 3$ we get
\begin{equation}
    p\left(\frac{\cos(2N\theta)}{3} \right) = -\frac{4}{9} q_N^2 \leq 0
\end{equation}
hence now the eigenvalue is at least $7/24$. Thus, the maximal eigenvalue of $Q_N$ is lower-bounded everywhere by $1/4$, as claimed. 
\vspace{0.1in}
\newline
To obtain a $Q_N$ in dimensions $n \geq 3$, we can simply lift the vector field $V_N$ to higher dimension. Consider the sphere $\mathbb{S}^m$ as a family of smaller spheres $\mathbb{S}^{m-1}_{\theta}$
\begin{equation}
    (\cos(\theta), \sin(\theta) x) \quad x \in \mathbb{S}^{m-1}, 0 \leq \theta \leq \pi
\end{equation}
with round metric $h_m = d\theta^2 + \sin^2(\theta) h_{m-1}$. On each of those, we can use $V_{N, m-1}$ for dimension $m-1$. At the poles $\theta = 0, \pi$, one adds cutoff $f, g$ supported near/away from the poles (similar to the $\chi$ we used in two dimensions), and then defines
\begin{equation}
    V_{N, m} = f_N(\theta) \partial_{\theta} + g_N(\theta) j_{\ast} V_{N, m-1}
\end{equation}
where $j(x) = (\cos(\theta), \sin(\theta) x)$ denotes the inclusion of $\mathbb{S}^{m-1}$ into $\mathbb{S}^m$. It can be checked exactly as above that this satisfies the required property, i.e. its norm goes to $0$ as $N$ grows, $Q_N$ is bounded by $1$ and finally its maximal eigenvalue $\lambda_{\max}$ is lower-bounded. Perhaps more interestingly, one can also construct perturbations $Q_N$ for which multiple eigenvalues are lower-bounded, which creates multiple directions of concavity. Based on the same construction as above and using ChatGPT Astra-6, the author managed to construct an example with $n-2$ uniformly positive eigenvalues on $\mathbb{S}^{n-1}$ satisfying the desired assumptions to preserve simplicity, which makes the 'concave subspace' of $Q_N$ codimension $2$ in the manifold. The construction is again iterative, and proceeds by projecting to the sphere an oscillating vector field $U$ which itself supplies $n-1$ positive eigenvalues (hence leaving $n-2$ on the sphere):
\begin{equation}
\begin{aligned}
    &V_N = \frac{1}{2nN}(I - x \otimes x) U_n(Nx) \\
    & U_n(y, z) = \left(U_{n-1}(y) - C_n \cos(\omega_n z) Z_{n-1}(y), \frac{\sin(2\omega_n z)}{2 \omega_n} \right) 
\end{aligned}
\end{equation}
where the constant $C_n$ is chosen small enough and the frequency $\omega_n$ a large enough integer (depending on $C_n$), while $Z_{n-1}$ denotes a smooth choice of eigenvector with the smallest eigenvalue $\lambda_{\min}$ of $sym \nabla U_{n-1}$. Then, the $n-2$ positive eigenvalues of $sym \nabla U_{n-1}$ in $sym \nabla U_{n}$ are preserved (using $C_n$ small), while on the two-dimensional subspace consisting of $Z_{n-1}(y)$ and $\partial_z$ the symmetric derivative is (up to an error $O(C_n)$ coupling that subspace with the other $n-2$ directions)
\begin{equation}
    sym \nabla U_{n} = \begin{pmatrix}
\lambda_{\min}(y)
  & \dfrac{C_n \omega_n}{2}\sin(\omega_n z) \\
\dfrac{C_n \omega_n}{2}\sin(\omega_n z)
  & \cos(2\omega_n z)
\end{pmatrix}.
\end{equation}
For large enough $\omega_n$, its largest eigenvalue is at least 1. Thus, together with the other $n-2$ positive eigenvalues (which are in the orthogonal complement of this subspace), this gives $n-1$ positive eigenvalues. Restricting the vector field to the sphere loses at most one eigenvalue, hence we get the desired $n-2$ dimensional subspace. Note that $V_N \rightarrow 0$ as $N \rightarrow \infty$ as desired. This tells us that the standard convexity method to invert the X-ray transform most likely fails for simple manifolds, even in higher dimensions. 

\section{Proof of the formula for the curvature tensor.}
Here we justify lemma \ref{lemma:CurvatureTensor} that we needed for the case $1 \leq \beta \leq b$ of our proof of simplicity. We recall the decomposition into unperturbed, oscillating and pointwise small parts:
\begin{equation}
\label{eq:CurvatureAppendix}
    \langle \mathcal{R}(J, \dot{\gamma}) \dot{\gamma}, J \rangle = |J|^2 - \frac{1}{2} \psi''\left( \frac{\pi/2 - r}{\epsilon} \right) Q_{x(t)}(X, X) \langle J, V \rangle^2 + \langle\mathcal{E}J, J\rangle \quad V \text{ some almost radial vector field}
\end{equation}
with $\mathcal{E} = O_b(\delta)$. The idea is quite simple. We already have a priori control on the perturbation 
\begin{equation}
    H(r) = \sin^2(r)h + \epsilon^2 \psi\left(\frac{\pi/2 - r}{\epsilon}\right) Q
\end{equation}
in the desired region. We recall that it satisfies $H(r) = h + O_b(\epsilon^2)$ while its radial derivative is $\partial H = O_b(\epsilon)$ and its tangential ones (to the sphere) are controlled by our estimates on $Q$. When computing the curvature tensor, we take up to two derivatives of the metric. The worst case scenario is when two radial derivatives land on $H$, in which case we end up with the dominating contribution
\begin{equation}
    \partial_r^2 \left( \epsilon^2 \psi\left(\frac{\pi/2 - r}{\epsilon}\right) Q_N \right) \sim \psi'' Q
\end{equation}
which corresponds to the oscillating contribution. The remaining contributions all have at most a single radial derivative, in which case they are bounded by $O(\epsilon)$ times the bound on the tangential derivatives of $Q$, which gives $O(\delta)$. This is basically the proof, but we give the formulas for completeness. 
\vspace{0.1in}
\newline
Radial component of the curvature: The Christoffel symbols are given in terms of $H$ by
\begin{equation}
    \Gamma^r_{ij} = \frac{-1}{2}(\partial_r H)_{ij} \quad \Gamma^i_{rj} = \frac{1}{2} H^{il}(\partial_r H)_{lj}
\end{equation}
Using $H(r)$ as above and tangential vector fields $X, Y$ to the sphere we get
\begin{equation}
\begin{aligned}
    \langle \mathcal{R}(X, \partial_r) \partial_r, Y \rangle &= \left( \frac{-1}{2} \partial_r^2 H + \frac{1}{4} (\partial_r H) H^{-1} (\partial_r H) \right)(X, Y) \\
    &=H(X, Y) - \frac{1}{2} \psi'' Q(X, Y) + O_b(\epsilon^2) \sqrt{h(X, X)} \sqrt{h(Y, Y)}
\end{aligned}
\end{equation}
\vspace{0.1in}
\newline
Mixed components of the curvature: In that case, use the Codazzi formula to relate it to the extrinsic curvature $\Pi = \frac{1}{2} \partial_r H$ of the spheres:
\begin{equation}
    \langle \mathcal{R}(X, Y) \partial_r, W \rangle = \nabla^H_X \Pi(Y, W) - \nabla^H_Y \Pi (X, W)
\end{equation}
By what we showed in step 2 (i.e. the bounds on $H$), the connection of $H$ and $h$ are close together
\begin{equation}
    |\nabla^H - \nabla^h| \leq C_b \epsilon \delta 
\end{equation}
while $\nabla^h \Pi = -\frac{\psi' \epsilon}{2} \nabla^h Q = O(\delta)$. Thus, we conclude that the mixed components are $O_b(\delta)$. 
\vspace{0.1in}
\newline
Tangential components of the curvature: Now we use the Gauss formula
\begin{equation}
    \langle \mathcal{R}(X, Y) V, W \rangle = H(\mathcal{R}^H(X, Y)V, W) - \Pi(Y, V) \Pi(X, W) + \Pi(X, V) \Pi(Y, W)
\end{equation}
Now again using the estimates on $H$ and its covariant derivative, we get that
\begin{equation}
    |\mathcal{R}^H - \mathcal{R}^h| \leq C_b \epsilon \delta
\end{equation}
Similarly the extrinsic curvature is $O(\epsilon)$, hence we get an overall error of $O_b(\epsilon \delta)$ from the constant curvature of $1$ of the sphere, which is the $|J|^2$ term above.
\vspace{0.1in}
\newline
This already gives the desired control of the curvature tensor. To get it to the form \ref{eq:CurvatureAppendix}, one simply chooses $V$ which is almost radial, i.e.
\begin{equation}
    |V - \partial_r| \leq C_b \epsilon
\end{equation}
over a bounded time along the geodesic, but which is also normal to the curve's tangent vector. By initially orthonormalizing $\partial_r$
\begin{equation}
    V(0) = \frac{\partial_r - \dot{r}(0) \dot{\gamma}(0)}{\sqrt{1 - \dot{r}(0)^2}}
\end{equation}
and then parallel transporting it along the curve (for finite time with $D_t \partial_r = O(\epsilon)$), we obtain such a vector field $V$. Then, applying the formulas above to $V$ versus $\partial_r$ only gives an error $O(\epsilon) \leq O(\delta)$. 

\label{App:ProofCurvature}
\end{document}